\documentclass[reqno,12pt]{amsart}
\usepackage{graphicx,epsfig}
\usepackage{amssymb,amsthm,amsmath}
\usepackage{latexsym}

\newtheorem{thm}{Theorem}
\newtheorem{cor}[thm]{Corollary}
\newtheorem{lem}{Lemma}

\input xy
\xyoption{all}

\def\N{\hbox{\font\dubl=msbm10 scaled 1200 {\dubl N}}}
\def\Q{\hbox{\font\dubl=msbm10 scaled 1200 {\dubl Q}}}
\def\R{\hbox{\font\dubl=msbm10 scaled 1200 {\dubl R}}}
\def\d{\,{\rm{d}}}

\title[Statistical properties of the Calkin--Wilf tree]{Statistical properties of the Calkin--Wilf tree: real an $p-$adic distribution}

\author[Giedrius Alkauskas, J\"orn Steuding]{Giedrius Alkauskas, J\"orn Steuding}

\begin{document}

\begin{abstract} We examine statistical properties of the Calkin--Wilf tree and give number-theoretical applications.
\end{abstract}

\maketitle

\section{A mean-value related to the Calkin--Wilf tree}

The Calkin--Wilf tree is generated by the iteration
$$
{a\over b}\quad\mapsto\quad {a\over a+b}\ ,\quad {a+b\over b},
$$
starting from the root ${1\over 1}$; the number ${a\over a+b}$ is called the left child of ${a\over b}$ and ${a+b\over b}$ the right child; we also say that ${a\over b}$ is the mother of its children. Recently, Calkin \& Wilf \cite{cw} have shown that this tree contains any positive rational number once and only once, each of which represented as a reduced fraction. The first iterations lead to
$$
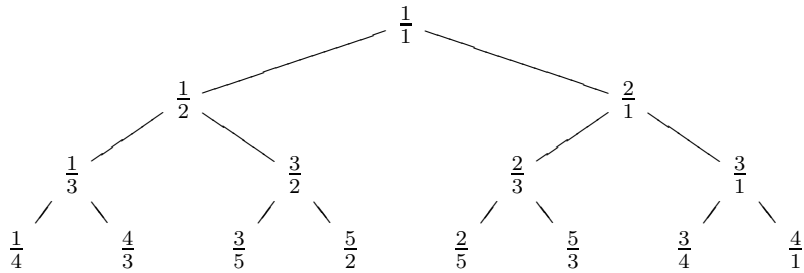
\xymatrix @R=.5pc @C=.5pc { & & & & & & & {1\over 1} & & & & & & & \\
& & & {1\over 2} \ar@{-}[urrrr] & & & & & & & & {2\over 1} \ar@{-}[ullll] & & & \\
& {1\over 3} \ar@{-}[urr] & & & & {3\over 2}\ar@{-}[ull] & & & & {2\over 3}\ar@{-}[urr] & & & & {3\over 1} \ar@{-}[ull] & \\
{1\over 4} \ar@{-}[ur] & & {4\over 3} \ar@{-}[ul] & & {3\over 5} \ar@{-}[ur] & & {5\over 2} \ar@{-}[ul] & & {2\over 5} \ar@{-}[ur]
& & {5\over 3} \ar@{-}[ul] & & {3\over 4} \ar@{-}[ur] & & {4\over 1} \ar@{-}[ul] }
$$
Reading the tree line by line, the Calkin--Wilf enumeration of $\Q^+$ starts with
$$
{1\over 1}\,,\ {1\over 2}\,,\ {2\over 1}\,,\ {1\over 3}\,,\
{3\over 2}\,,\ {2\over 3}\,,\ {3\over 1}\,,\ {1\over 4}\,,\
{4\over 3}\,,\ {3\over 5}\,,\ {5\over 2}\,,\ {2\over 5}\,,\
{5\over 3}\,,\ {3\over 4}\,,\ {4\over 1}\,,\ \ldots
$$
As recently pointed out by Reznick \cite{reznick}, this sequence was already investigated by Stern \cite{stern} in 1858. This sequence satisfies also the iteration
$$
x_1=1,\quad x_{n+1}=1/(2[x_n]+1-x_n),
$$
where $[x]$ denotes the largest integer $\leq x$; this observation is due to Newman (cf. \cite{new}), answering a question of D.E. Knuth, resp. Vandervelde \& Zagier (cf. \cite{sloane}).
\par

The Calkin--Wilf enumeration of the positive rationals has many interesting features. For instance, it encodes the hyperbinary representations of all positive integers (see \cite{cw}). Furthermore, it can be used as model for the game {\it Euclid} first formulated by Cole \& Davie \cite{euclid}; see Hofmann, Schuster \& Steuding \cite{hss}. In this short note we are concerned with statistical properties of the Calkin--Wilf tree.
\par

We write the $n$th generation of the Calkin--Wilf tree as ${\mathcal C}{\mathcal W}^{(n)}=\{x_j^{(n)}\}_j$, where the $x_n^{(j)}$ are the elements ordered according to their appearance in the $n$th line of the Calkin--Wilf tree. So $\Q^+=\bigcup_{n=1}^\infty{\mathcal C}{\mathcal W}^{(n)}$. Obviously, ${\mathcal C}{\mathcal W}^{(n)}$ consists of $2^{n-1}$ elements. Denote by $\Sigma(n)$ the sum of all elements of the $n$th generation of the Calkin--Wilf tree,
$$
\Sigma(n)=\sum_{j=1}^{2^{n-1}}x_j^{(n)}.
$$
Our first result gives the mean-value of the elements of the $n$th generation of the Calkin--Wilf tree:

\begin{thm}\label{uno}
For any $n\in\N$,
$$
\Sigma(n)=3\cdot 2^{n-2}-{1\over 2}.
$$
\end{thm}

\noindent This result may be interpreted as follows. We observe that $x_1^{(n)}={1\over n}$ and $x_{2^{n-1}}^{(n)}={n\over 1}$ for all $n\in\N$, and thus ${\mathcal C}{\mathcal W}^{(n)}$ is supported on an unbounded set as $n\to\infty$. However, the average value of the $2^{n-1}$ elements of the $n$th generation ${\mathcal C}{\mathcal W}^{(n)}$ is approximately ${3\over 2}$, which is, surprisingly, a finite number. This has a simple explanation: in some sense, {\it small} values are taken in earlier generations than {\it large} values. For instance, in each generation ${\mathcal C}{\mathcal W}^{(n)}$ takes as many values form the interval $(0,1)$ as from $(1,\infty)$. This result was also recently proved by Reznick \cite{reznick}; his proof differs slightly from our argument.\footnote{The problem of determining the average value of the Calkin--Wilf tree was posed by the second named author as a problem in the problem session of the IV International conference on analytic and probabilistic number theory in Palanga 2006; an independent solution was given by Eduard Wirsing.}
\par\bigskip

\noindent {\bf Proof} by induction on $n$. The statement of the theorem is correct for $n=1$ and $n=2$. Now suppose that $n\geq 3$. In order to prove the statement for $n$ we first observe a symmetry in the Calkin--Wilf tree with respect to its middle: for $n\geq 2$,
\begin{equation}\label{sym}
x_j^{(n)}={a\over b}\qquad\iff\qquad x_{2^{n-1}+1-j}^{(n)}={b\over
a}\ ;
\end{equation}
this is easily proved by another induction on $n$ (and we leave its simple verification to the reader). Further, we note that $x_j^{(n)}\leq 1$ if and only if $j$ is odd; here equality holds if and only if $n=1$.
\par

Now we start to evaluate $\Sigma(n)$. For this purpose we compute
$$
y_j^{(n)}:=\left\{\begin{array}{c@{\quad}l} x_j^{(n)}+ x_{2^{n-1}-j}^{(n)} & \mbox{for}\
j=1,2,\ldots,2^{n-2}-1,\\ x_j^{(n)}+ x_{2j}^{(n)}  &
\mbox{for}\ j=2^{n-2},\end{array}\right.
$$
and add these values over $j=1,2,\ldots,2^{n-2}$. Clearly, $\Sigma(n)=\sum_{j=1}^{2^{n-2}}y_j^{(n)}$.
\par

First, assume that $j$ is odd. Then both, $x_j^{(n)}$ and $x_{2^{n-1}-j}^{(n)}$ are strictly less than $1$. In view of (\ref{sym}) the mothers of $x_j^{(n)}$ and $x_{2^{n-1}-j}^{(n)}$ are of the form ${a\over b}$ and ${b\over a}$, respectively. Hence,
$$
x_j^{(n)}+ x_{2^{n-1}-j}^{(n)}={a\over a+b}+{b\over a+b}
$$
and thus we find $y_j^{(n)}=1$ in this case.
\par

Next, we consider the case that $j$ is even. Then both, $x_j^{(n)}$ and $x_{2^{n-1}-j}^{(n)}$ are strictly greater than $1$. If the mothers of $x_j^{(n)}$ and $x_{2^{n-1}-j}^{(n)}$ are of the form ${a\over b}$ and ${a'\over b'}$, respectively, then
$$
x_j^{(n)}={a+b\over b}=1+{a\over b}\qquad\mbox{and}\qquad
x_{2^{n-1}-j}^{(n)}=1+{a'\over b'}.
$$
Hence, we find for their sum
$$
x_j^{(n)}+ x_{2^{n-1}-j}^{(n)}=2+{a\over b}+{a'\over b'}
$$
and so $y_j^{(n)}=2+y_k^{(n-1)}$, where $y_k^{(n-1)}$ is either the sum of two elements $x_k^{(n-1)}$ and $x_{2^{n-2}-k}^{(n-1)}$ or the sum of $x_{2^{n-3}}^{(n-1)}$ and $x_{2^{n-2}}^{(n-1)}$.
\par

It remains to combine both evaluations. Since both cases appear equally often, namely each $2^{n-3}$ times, we obtain the recurrence formula
$$
\Sigma(n)=\Sigma(n-1)+(1+2)\cdot 2^{n-3},
$$
being valid for $n\geq 3$. This implies the assertion of the theorem. $\blacksquare$

\section{An application to finite continued fractions}

Theorem \ref{uno} has a nice number-theoretical interpretation. It is well-known that each positive rational number $x$ has a representation as a finite (regular) continued fraction
$$
x=a_0+{1\over a_1+{ \atop \displaystyle{\ddots { \atop +{1 \over {a_{m-1}+{1\over a_m}}}}}}}
$$
with $a_0\in\N\cup\{0\}$ and $a_j\in\N$ for some $m\in\N\cup\{0\}$. In order to have a unique representation, we assume that $a_m\geq 2$ if $m\in\N$. We shall use the standard notation $x=[a_0,a_1,\ldots,a_m]$. Continued fractions are of special interest in the theory of diophantine approximation.
\par

As Bird, Gibbons \& Lester \cite{bgl} showed, the $n$th generation of the Calkin--Wilf tree consists exactly of those rationals having a continued fraction expansion $[a_0,a_1,\ldots,a_m]$ for which the sum of the partial quotients $a_j$ is constant $n$, the continued fractions of even length in the left subtree, and the continued fractions with odd length in the right subtree. Thus Theorem \ref{uno} yields

\begin{cor}\label{appl}
For any $n\in\N$,
$$
2^{1-n}\sum_{a_0+a_1+\ldots+a_m=n}[a_0,a_1,\ldots,a_m]={3\over
2}-2^{-n}.
$$
\end{cor}

\noindent One can use the approach via continued fractions to locate any positive rational in the tree. This observation is due to Bird, Gibbons \& Lester \cite{bgl} (actually, their reasoning is based on Graham, Knuth \& Patasnik \cite{graham} who gave such a description for the related Stern--Brocot tree). Given a reduced fraction $x$ in the Calkin--Wilf tree with continued fraction expansion
$$
x=[a_0,a_1,\ldots,a_{m-2},a_{m-1},a_m],
$$
we associate the path
\begin{eqnarray*}
{\sf L}^{a_m-1}{\sf R}^{a_{m-1}}{\sf L}^{a_{m-2}}\cdots{\sf L}^{a_1}{\sf R}^{a_0}&\mbox{if}& m\ \mbox{is odd, and}\\
{\sf R}^{a_m-1}{\sf L}^{a_{m-1}}{\sf R}^{a_{m-2}}\cdots{\sf R}^{a_1}{\sf L}^{a_0}&\mbox{if}& m\ \mbox{is even};
\end{eqnarray*}
note that $a_m-1\geq 1$ for $m\in\N$. The notation ${\sf R}^a$ with $a\in\N\cup\{0\}$ means: $a$ steps to the right, whereas ${\sf L}^b$ with $b\in\N\cup\{0\}$ stands for $b$ steps to the left. Then, starting from the root ${1\over 1}$ and following this path from left to right, we end up with the element $x$. This follows easily from the iteration with which the tree was build; notice that this claim is essentially already contained in Lehmer \cite{lehmer} (this was also observed by Reznick \cite{reznick}).

\begin{cor}\label{korol}
Given any non-empty interval $(\alpha,\beta)$ in $\R^+$, and any finite path in the Calkin--Wilf tree, there exists a continuation
of this path which contains a rational number from the interval $(\alpha,\beta)$.
\end{cor}

\noindent {\bf Proof.} We expand $\alpha$ and $\beta$ into continued fractions, $\alpha=[a_0,a_1,\ldots]$ and $\beta=[b_0,b_1,\ldots]$, say. Let $k$ be the least index such that $a_k\neq b_k$. According to the parity of $k$ we have $a_k<b_k$ (if $k$ is even) or $a_k>b_k$ (if $k$ is odd). Without loss of generality we may assume that $\vert b_k-a_k\vert\geq 2$ (since otherwise we may consider a subinterval of $(\alpha,\beta)$). Moreover we may suppose that the path in question is starting from the root and is given in the form ${\sf L}^{c_m-1}{\sf R}^{c_{m-1}}{\sf L}^{c_{m-2}}\cdots{\sf L}^{c_1}{\sf R}^{c_0}$ (the other case may be treated analogously). Then we construct a rational number $x$ by assigning the finite continued fraction
$$
x=[a_0,a_1,\ldots,a_{k-1},x_k,x_{k+1},c_0,c_1,\ldots,c_{m-2},c_{m-1},c_m],
$$
where $x_k:=\min\{a_k,b_k\}+1$ and $x_{k+1}$ denotes the string $1$ if $k$ is odd, resp. $1,1$ if $k$ is even. Since $b_j=a_j$ for $0\leq j<k$ and
$$
\min\{a_k,b_k\}<x_k<\max\{a_k,b_k\},
$$
it follows that $\alpha<x<\beta$. Since the length of the continued fraction expansion has the same parity as $m$ (thanks to the definition of $x_{k+1}$), the element $x$ can be reached by the path ${\sf L}^{c_m-1}{\sf R}^{c_{m-1}}{\sf L}^{c_{m-2}}\cdots{\sf L}^{c_1}{\sf R}^{c_0}$. This proves the corollary. $\blacksquare$

\section{A random walk on the Calkin--Wilf tree}

Starting with $X_1={1\over 1}$, we define a sequence of random variables
by the following iteration: if $X_n={a\over b}$, then $X_{n+1}={a\over a+b}$ with probability ${1\over 2}$ and
$X_{n+1}={a+b\over b}$ with probability ${1\over 2}$. The sequence $\{X_n\}$ may be regarded as a random walk on the
Calkin--Wilf tree where $n$ is a discrete time parameter.

\begin{thm}
Let $(\alpha,\beta)$ be any non-empty interval in $\R^+$. Then, with probabilty $1$, the random walk $\{X_n\}$
visits the interval $(\alpha,\beta)$, i.e., with probability $1$, there exists $m\in\N$ such that $x_m\in(\alpha,\beta)$.
\end{thm}

\noindent {\bf Proof.} The interval $(\alpha,\beta)$ contains a non-empty subinterval $[A,B]$ such that for any $\zeta\in[A,B]$ the initial partial quotients $c_0,c_1,\ldots,c_m$ are identical: $\zeta=[c_0,c_1,\ldots,c_m,\ldots]$. Hence, with the interval $[A,B]$ we may associate a path pattern ${\sf L}^{c_m-1}{\sf R}^{c_{m-1}}{\sf L}^{c_{m-2}}\cdots{\sf L}^{c_1}{\sf R}^{c_0}$ in the Calkin-Wilf tree such that any path in the tree starting from the root and ending with ${\sf L}^{c_m-1}{\sf R}^{c_{m-1}}{\sf L}^{c_{m-2}}\cdots{\sf L}^{c_1}{\sf R}^{c_0}$ points to an element in $[A,B]$. Since the probability is ${1\over 2}$ for both ${a\over b}\mapsto {a\over a+b}$ and ${a\over b}\mapsto {a+b\over b}$, each pattern of fixed length $m$ appears with the same probability and so we may restrict on the path pattern $R^k$.
\par

In the case $k=1$ we find in each generation exactly one which ends with $R$ but does not contain any $R$ before (actually, this is $L^{n-1}R$ in generation $n$). Adding up all probabilities for these paths, we get
$$
\sum_{j=1}^\infty\left({1\over 2}\right)^j=1,
$$
and so, with probability $1$, the random walk $X_n$ will go to the right child for some $n$. Now assume that the statement is true for $k$. We shall show that then it is also true for $k+1$. For each path of the form $XR^k$, where $X$ is any combination of powers of $L$ and $R$, there are two paths $XR^kL$ and $XR^{k+1}$, so by induction the probability that the random walk eventually follows the path $R^{k+1}$ is at least ${1\over 2}$. However, for each path $XR^k$ one also has to consider the subtrees starting from $XR^kL^d$ for $d=1,2,\ldots$, each of which containing paths which end $R^{k+1}$. By self-similarity, the probability that the random walk eventually follows the path $R^{k+1}$ is
$$
\sum_{j=1}^\infty\left({1\over 2}\right)^j\sum_{d=1}^\infty \left({1\over 2}\right)^d=1.
$$
This proves the theorem. $\blacksquare$

\section{Statistical properties of the Calkin--Wilf tree}

In view of Corollary \ref{appl} it is interesting to have a better understanding of the statistics of the Calkin--Wilf tree. The following theorem gives the limit distribution function in explicit form.

\begin{thm}\label{wno}
Let $F_{n}(x)$ denote the distribution function of the $n-$th
generation, i.e.,
$$
F_{n}(x)=2^{1-n}\#\{j\,:\,x_{j}^{(n)}\leq x\}.
$$
Then uniformly $F_{n}(x)\rightarrow F(x)$, where
$$
F([a_{0},a_{1},a_{2},a_{3},...])=1-2^{-a_{0}}+2^{-(a_{0}+a_{1})}-2^{-(a_{0}+a_{1}+a_{2})}+...
$$
(for rational numbers $x=[a_0,a_1,\ldots]$ this series terminates at the last non-zero partial quotient of the continued fraction). Thus, $F(0)=0$, $F(\infty)=1$, and $F(x)$ is a monotonically increasing function. Moreover, $F(x)$ is continuous and singular, i.e., $F'(x)=0$ almost everywhere.
\end{thm}

\noindent {\bf Proof.} Let $x\geq1$. One half of the fractions in the $n+1-$st generation do not exceed $1$, and hence also do not exceed $x$. Further,
$$
\frac{a+b}{b}\leq x\quad \iff\quad \frac{a}{b}\leq x-1.
$$
Hence,
$$
2F_{n+1}(x)=F_{n}(x-1)+1,\quad n\geq 1.
$$
Now assume $0<x<1$. Then
$$
\frac{a}{a+b}\leq x\quad\iff\quad \frac{a}{b}\leq\frac{x}{1-x}.
$$
Therefore,
$$
2F_{n+1}(x)=F_{n}\Big{(}\frac{x}{1-x}\Big{)}.
$$
The distribution function $F$, defined in the formulation of the theorem, satisfies the functional equation
$$
2F(x)=\left\{\begin{array}{c@{\qquad}l} F(x-1)+1 & \mbox{if}\quad x\geq 1,
\\ F({x\over 1-x}) & \mbox{if}\quad 0<x<1. \end{array}\right.
$$
For instance, the second identity is equivalent to $2F(\frac{t}{t+1})=F(t)$ for all positive $t$. If $t=[b_{0},b_{1},...]$, then $\frac{t}{t+1}=[0,1,b_{0},b_{1},..]$ for $t\geq 1$, and $\frac{t}{t+1}=[0,b_{1}+1,b_{2},...]$ for $t<1$, and the statement follows immediately.
\par

Now define $\delta_{n}(x)=F(x)-F_{n}(x)$. In order to prove the first assertion of the theorem, the uniform convergence $F_n\to F$, it is sufficient to show that
\begin{equation}\label{vienas}
\sup_{x\geq 0}|\delta_{n}(x)|\leq 2^{-n}.
\end{equation}
It is easy to see that the assertion is true for $n=1$. Now suppose the estimate is true for $n$. In view of the functional equation for both $F_n(x)$ and $F(x)$, we have
$$
2\delta_{n+1}(x)=\delta_{n}\Big{(}\frac{x}{1-x}\Big{)}
$$
for $0<x<1$, which gives $\sup_{0\leq x<1}|\delta_{n+1}(x)|\leq 2^{-n-1}$. Moreover, we have
$$
2\delta_{n+1}(x)=\delta_{n}(x-1)
$$
for $x\geq1$, which yields the same bound for $\delta_n(x)$ in the range $x\geq 1$. This proves (\ref{vienas}).
\par

Clearly, $F$, as a distribution function, is monotonic; obviously, it is also continuous. It remains to prove that $F(x)$ is singular. Given an irrational number $\alpha=[a_{0},a_{1},a_{2},...]$, we consider the sequence
$$
\alpha_{n}=[a_{0},a_{1},...,a_{n-1},a_{n}+1,a_{n+1},...];
$$
obviously, $\alpha_n$ is the real number which is defined by the continued fraction expansion of $\alpha$, where the $n$th partial quotient $a_n$ is replaced by $a_n+1$. Depending on the parity of $n$, $\alpha_{n}$ is less than or greater than $\alpha$. Thus, any real number $y$, which is sufficiently close to $\alpha$, is contained between two terms of the sequence, $\alpha_{L}$ and $\alpha_{L+2}$ say. Then
$$
\Big{|}\frac{F(y)-F(\alpha)}{y-\alpha}\Big{|}\leq
\Big{|}\frac{F(\alpha_{L})-F(\alpha)}{\alpha_{L+2}-\alpha}\Big{|}.
$$
From the explicit form of $F$ we deduce
$$
|F(\alpha_{L})-F(\alpha)|\leq\frac{1}{2}2^{-(a_{0}+a_{1}+...+a_{L})}.
$$
On the other hand,
\begin{eqnarray*}
|\alpha_{L+2}-\alpha|&\geq& ([a_{1},a_{2},...,a_{L+2}+1,...]-[a_{1},a_{2},...,a_{L+2},...])(a_{0}+1)^{-2}\\
&\geq & \Big{(}(a_{0}+1)(a_{1}+1)...(a_{L+2}+1)\Big{)}^{-2}
\end{eqnarray*}
by induction. Thus,
$$
\Big{|}\frac{F(y)-F(\alpha)}{y-\alpha}\Big{|}
\leq 2^{1-(a_{0}+a_{1}+...+a_{L})}\prod_{i=1}^{L+2}(a_{i}+1)^{2}.
$$
The theorem of Khinchin (\cite{khin}, p. 86, implies that $\prod_{i=1}^{n}(a_{i}+1)^{1/n}$ tends to a fixed constant limit almost everywhere. On the other hand, the same reasoning shows that ${1\over n}\sum_{i=1}^{n}a_{n}$ tends to infinity for almost all $x$. Thus, almost everywhere the limit
$$
\lim_{y\rightarrow\alpha}(F(y)-F(\alpha))(y-\alpha)^{-1}
$$
exists and is equal $0$. This finishes the proof of the theorem. $\blacksquare$
\par\medskip

By the same argument as for the singular behaviour of $F$ we can show that $F'(\frac{\sqrt{5}+1}{2})=\infty$. Actually, the terms of ${\mathcal CW}^{(n)}$ are densely concentrated around numbers with $F'(x)=\infty$ and scarcely around those where $F'(x)=0$. The value of $F(x)$ is rational iff $x$ is either rational or quadratic irrationality, e.g.
$$
F(1)=\frac{1}{2},\quad F(\sqrt{2})=\frac{3}{5},\quad F((\sqrt{5}+1)/2)=\frac{2}{3}.
$$
This follows immediately from Lagrange's theorem which characterizes the quadratic irrationals by their eventually periodic continued fraction expansion.
For Euler's number $e=[2,\overline{1,2n,1}]$ we find that $F(e)$ can be expressed in terms of special values of Jacobi theta functions.

\section{Characteristics of the distribution function}

In view of Corollary $\ref{appl}$, the mean of the distribution function $F$ is ${3\over 2}$. Since $F$ has a tail of exponential decay, more precisely $1-F(x)=O(2^{-x})$, it follows that all moments exist. For $k\in\N_0$, the $k$th moment is defined by
$$
M_{k}=\int\limits_{0}^{\infty}x^{k}\d F(x).
$$
In order to give an asymptotic formula for $M_k$ let
$$
m_{k}=\int\limits_{0}^{\infty}\Big{(}\frac{x}{x+1}\Big{)}^{k}\d F(x)
$$
We will see that the generating function of $m_{k}$ has some interesting properties. Let $\omega(x)$ be a continuous function of at most power growth: $\omega(x)\ll x^{T}$ as $x\to\infty$. By the functional equation for $F$ we find $F(x+n)=1-2^{-n}+2^{-n}F(x)$, $x\geq0$. Hence
\begin{eqnarray*}
\int\limits_{0}^{\infty}\omega(x)\d F(x)&=&
\sum\limits_{n=0}^{\infty}\int\limits_{0}^{1}\omega(x+n)\d F(x+n)\\
&=& \int\limits_{0}^{1}\sum\limits_{n=0}^{\infty}\frac{\omega(x+n)}{2^{n}}\d F(x);
\end{eqnarray*}
these integrals exist in view of our assumptions and the fact that $F(x)$ has a tail of exponential decay. Let $x=\frac{t}{t+1}$ for $t\geq0$. Since $F(\frac{t}{t+1})=\frac{1}{2}F(t)$, this change of variables gives
$$
\int\limits_{0}^{\infty}\omega(x)\d F(x)=
\sum\limits_{n=0}^{\infty}
\int\limits_{0}^{\infty}\frac{\omega(\frac{t}{t+1}+n)}{2^{n+1}}\d F(t)
$$
(All changes of order of summation and integration are justified by the condition we put on $\omega(x)$). Now let $\omega(x)=x^{L}$ for some $L\in\mathbb{N}_{0}$ and define
$$
b_{s}=\sum_{n=0}^{\infty}\frac{n^{s}}{2^{n+1}}.
$$
Then
$$
\int\limits_{0}^{\infty}x^{L}\d F(x)=\int\limits_{0}^{\infty}\sum\limits_{i=0}^{L}\Big{(}\frac{x}{x+1}\Big{)}^{i}\binom{L}{i}b_{L-i}\d F(x),
$$
whence the relation
\begin{equation}\label{Mm}
M_{L}=\sum\limits_{i=0}^{L}m_{i}\binom{L}{i}b_{L-i}
\end{equation}
for $L\in\N_0$. The generating function of the sequence of the $b_s$ is given by
$$
b(t)=\sum\limits_{L=0}^{\infty}\frac{b_{L}}{L!}t^{L}=\sum\limits_{L=0}^{\infty}\sum\limits_{n=0}^{\infty}\frac{n^{L}t^{L}}{2^{N+1}L!}=\sum\limits_{n=0}^{\infty}\frac{e^{nt}}{2^{n+1}}=\frac{1}{2-e^{t}}.
$$
Denote by $M(t)$ and $m(t)$ the corresponding generating functions of the coefficients $M_{k}$ and $m_{k}$, respectively. Then we can rewrite (\ref{Mm}) as
$$
M(t)=\sum\limits_{L=0}^{\infty}\frac{M_{L}}{L!}t^{L}=\frac{1}{2-e^{t}}
\sum\limits_{L=0}^{\infty}\frac{m_{L}}{L!}t^{L}=\frac{1}{2-e^{t}}m(t).
$$
The function $m(t)$ is entire, and $M(t)$ has a positive radius of convergence. This already allows us to find approximate values of the moments $M_{L}$.

\begin{thm} For $L\in\N_0$,
\begin{eqnarray*}
M_{L}&=&\frac{m(\log 2)}{2\log 2}\Big{(}\frac{1}
{\log2}\Big{)}^{L}L!+O_{\varepsilon}\Big{(}((4\pi^{2}+(\log 2)^{1/2}-\varepsilon)^{-L}\Big{)}L!\\
&=&
\Big{(}\frac{m(\log 2)}{2\log 2}\Big{(}\frac{1}
{\log2}\Big{)}^{L}+O(6.3^{-L})\Big{)}L!
\end{eqnarray*}
\end{thm}

\noindent {\bf Proof.} By Cauchy's formula, for any sufficiently small $r$,
$$
M_L={L!\over 2\pi i}\int_{\vert z\vert=r}{M(z)\over z^{L+1}}\d z.
$$
Changing the path of integration, we get by the calculus of residues
$$
M_L=-\mbox{Res}_{z=\log 2}\left({m(z)\over (2-e^z)z^{L+1}}\right)-{L!\over 2\pi i}\int_{\vert z\vert=R}{m(z)\over 2-e^z}{\d z\over z^{L+1}},
$$
where $R$ satisfies $\log 2<R<\vert \log 2+2\pi i\vert$ (which means that there is exactly one simple pole of the integrand located in the interior of the circle $\vert z\vert=R$). It is easily seen that the residue coincides with the main term in the formula of the lemma; the error term follows from estimating the integral. $\blacksquare$
\par\bigskip

We obtain the inverse to the linear equations (\ref{Mm}):
$$
m_{L}=M_{L}-\sum\limits_{s=0}^{L-1}M_{s}\binom{L}{s}
$$
for $L\in\N_0$. Since $b(t)(2-e^{t})=1$, the coefficients $b_{L}$ can be
calculated recursively
$$
b_{L}=\sum_{s=0}^{L-1}\binom{L}{s}b_{s}.
$$
Thus, $b_{0}=1,b_{1}=1,b_{2}=3,b_{3}=13,b_{4}=75,b_{5}=541.$
\par

We proceed with a property of the function $m(t)$ which reflects
the symmetry of the distribution function: $F(y)+F(1/y)=1$.
Unfortunately, this property is still insufficient for determining
the coefficients $m_{L}$. As a matter of fact,
\begin{eqnarray*}
m_{L}=\int\limits_{0}^{\infty}\Big{(}\frac{x}{x+1}\Big{)}^{L}\d F(x)&=&-\int\limits_{0}^{\infty}\Big{(}\frac{1/x}{1/x+1}\Big{)}^{L}\d F(1/x)\\
&=&\int\limits_{0}^{\infty}\Big{(}\frac{1}{x+1}\Big{)}^{L}\d F(x).
\end{eqnarray*}
Since
$$
\Big{(}\frac{x}{x+1}\Big{)}^{L}=\Big{(}\frac{x+1-1}{x+1}\Big{)}^{L}=
\sum\limits_{s=0}^{L}\binom{L}{s}(-1)^{L-s}\Big{(}\frac{1}{x+1}\Big{)}^{L-s},
$$
this gives
$$
m_{L}=\sum\limits_{s=0}^{L}\binom{L}{s}(-1)^{s}m_{s}
$$
for $L\geq0$. For example, $m_{1}=m_{0}-m_{1}$, which gives $m_{1}={1\over 2}$ (since $m_0=1$), and thus $M_{1}={3\over 2}$ (see Theorem \ref{uno}). For the other coefficients we only get linear relations. Thus, $2m_{3}=-{1\over 2}+3m_{2}$. In terms of $m(t)$ the recursion formula above yields the identity
$$
m(t)=m(-t)e^{t}.
$$
We conclude this chapter with the result, which uniquely determines the function $m(t)$ (along with the condition $m(0)=1$).

\begin{thm}
The function $m(s)$ satisfies the integral equation
\begin{eqnarray*}
m(-s)=(2e^{s}-1)\int\limits_{0}^{\infty}m'(-t)J_{0}(2\sqrt{st})\d t,
\quad s\in\mathbb{R}_{+},
\end{eqnarray*}
where $J_{0}(*)$ stands for the Bessel function:
$$
J_{0}(z)=\frac{1}{\pi}\int\limits_{0}^{\pi}\cos(z\sin x)\d x
$$
\end{thm}

\noindent The proof of this theorem and the solution of this integral equation, and thus the explicit description of the moments will be given in a subsequent paper.

\section{$p-$adic distribution}

In the previous sections, we were interested in the distribution of the $n$th generation of the tree $\mathcal{C}\mathcal{W}$ in the field of real numbers. Since the set of non-equivalent valuations of $\mathbb{Q}$ contains a valuation, associated with any prime number $p$, it is natural to consider the distribution of the set of each generation in the field of $p-$adic numbers $\mathbb{Q}_{p}$. In this case we have an ultrametric inequality, which implies that two circles are either co-centric, or do not intersect. We define
\begin{eqnarray*}
F_{n}(z,\nu)=2^{-n+1}\#\{\frac{a}{b}\in\mathcal{C}\mathcal{W}^{(n)}:\text{ord}_{p}(\frac{a}{b}-z)\geq\nu\},\quad
z\in\mathbb{Q}_{p},\quad\nu\in\mathbb{Z}.
\end{eqnarray*}
(When $p$ is fixed, the subscript $p$ in $F_{n}$ is omitted). Note that in order to calculate $F_{n}(z,\nu)$ we can confine to the case $\text{ord}_{p}(z)<\nu$; otherwise $\text{ord}_{p}(\frac{a}{b}-z)\geq\nu\Leftrightarrow\text{ord}_{p}(\frac{a}{b})\geq\nu$. We shall calculate the limit distribution
$\mu_{p}(z,\nu)=\lim_{n\rightarrow\infty}F_{n}(z,\nu)$, and also some characteristics of it, e.g. the zeta function
$$
Z_{p}(s)=\int\limits_{u\in\mathbb{Q}_{p}}|u|^{s}d\mu_{p},\quad s\in\mathbb{C},\quad z\in\mathbb{Q}_{p},
$$
where $|*|$ stands for the $p-$adic valuation.
\par

To illustrate how the method works, we will calculate the value of $F_{n}$ in two special cases. Let $p=2$, and let $E(n)$ be the number of rational numbers in the $n$th generation with one of $a$ or $b$ being even, and let $O(n)$ be the corresponding number fractions with both $a$ and $b$ odd. Then $E(n)+O(n)=2^{n-1}$. Since $\frac{a}{b}$ in the $n$th generation generates $\frac{a}{a+b}$ and $\frac{a+b}{b}$ in the $(n+1)$st generation, each fraction $\frac{a}{b}$ with one of the $a$, $b$ even will generate one fraction with both numerator and denominator odd. If both $a$, $b$ are odd, then their two offsprings will not be of this kind. Therefore, $O(n+1)=E(n)$. Similarly, $E(n+1)=E(n)+2O(n)$. This gives the recurrence $E(n+1)=E(n)+2E(n-1)$, $n\geq 2$, and this implies
$$
E(n)=\frac{2^{n}+2(-1)^{n}}{3},\quad O(n)=\frac{2^{n-1}+2(-1)^{n-1}}{3}, \quad \mu_{2}(0,0)=\frac{2}{3}.
$$
(For the last equality note that $\frac{a}{b}$ and $\frac{b}{a}$ simultaneously belong to $\mathcal{C}\mathcal{W}^{(n)}$, and so the number of fractions with $\text{ord}_{2}(*)>0$ is $E(n)/2$). We will generalize this example to odd prime $p\geq 3$. Let $L_{i}(n)$ be the part of the fractions in the $n$th generations such that $ab^{-1}\equiv i\bmod\,p$ for $0\leq i\leq p-1$ or $i=\infty$ (that is, $b\equiv0\bmod\,p$). Thus,
$$
\sum_{i\in\mathbb{F}_{p}\cup\infty}L_{i}(n)=1;
$$
in other words, $L_{i}(n)=F_{n}(i,1)$. For our later investigations we need a result from the theory of finite Markov chains.

\begin{lem}
Let $\textbf{A}$ be a matrix of a finite Markov chain with $s$
stages. That is, $a_{i,j}\geq0$, and $\sum_{j=1}^{s}a_{i,j}=1$ for
all $i$. Suppose that $\mathbf{A}$ is irreducible (for all pairs
$(i,j)$, and some $m$, the entry $a_{i,j}^{(m)}$ of the matrix
$\mathbf{A}^{m}$ is strictly positive), acyclic and recurrent
(this is satisfied, if all entries of $\mathbf{A}^{m}$ are
strictly positive for some $m$). Then the eigenvalue $1$ is
simple, if $\lambda$ is another eigenvalue, then $|\lambda|<1$,
and $\mathbf{A}^{m}$, as $m\rightarrow\infty$, tends to the matrix
$\mathbf{B}$, with entries $b_{i,j}=\pi_{j}$, where
$(\pi_{1},...,\pi_{s})$ is a unique left eigenvector with
eigenvalue $1$, such that $\sum_{j=1}^{s}\pi_{j}=1$.
\end{lem}

\noindent A proof of this lemma can be found in \cite{karl},
Section 3.1., Theorem 1.3.

\begin{thm}
$\mu_{p}(z,1)=\frac{1}{p+1}$ for $z\in\mathbb{Z}_{p}$.
\end{thm}

\noindent {\bf Proof. } Similarly as in the above example, a fraction $\frac{a}{b}$ from the $n$th generation generates $\frac{a}{a+b}$ and
$\frac{a+b}{b}$ in the $(n+1)$st generation, and it is routine to check that
\begin{eqnarray}
L_{i}(n+1)=\frac{1}{2}L_{\frac{i}{1-i}}(n)+\frac{1}{2}L_{i-1}(n)\quad\mbox{for}\quad
i\in\mathbb{F}_{p}\cup\{\infty\}\label{padic},
\end{eqnarray}
(Here we make a natural convention for $\frac{i}{1-i}$ and $i-1$, if $i=1$ or $\infty$). In this equation, it can happen that $i-1\equiv\frac{i}{1-i}\bmod\,p$; thus, $(2i-1)^{2}\equiv -3\bmod\,p$. The recurrence for this particular $i$ is to be understood in the obvious way, $L_{i}(n+1)=L_{i-1}(n)$.  Therefore, if we denote the vector-column $(L_{\infty}(n),L_{0}(n),...,L_{p-1}(n))^{T}$ by $\mathbf{v}_{n}$, and if $\mathcal{A}$ is a matrix of the system $(\ref{padic})$, then $\mathbf{v}_{n+1}=\mathcal{A}\mathbf{v}_{n}$, and hence
$$
\mathbf{v}_{n}=\mathcal{A}^{n-1}\mathbf{v}_{1},
$$
where $\mathbf{v}_{1}=(0,0,1,0,...,0)^{T}$. In any particular case, this allows us two find the values of $L_{i}$ explicitly. For example, if $p=7$, the characteristic polynomial is
$$
f(x)=\frac{1}{16}(x-1)(2x-1)(2x^{2}+1)(4x^{4}+2x^{3}+2x+1).
$$
The list of roots is
$$
\alpha_{1}=1,\quad\alpha=\frac{1}{2},\quad\alpha_{3,4}=\pm\frac{i}{\sqrt{2}},\quad\alpha_{5,6,7,8}=\frac{-1-\sqrt{17}}{8}\pm\frac{\sqrt{1+\sqrt{17}}}{2\sqrt{2}},
$$
(with respect to the two values for the root $\sqrt{17}$), the matrix is diagonalisible, and the Jordan normal form gives the expression
$$
L_{i}(n)=\sum_{s=1}^{8}C_{i,s}\alpha_{s}^{n}.
$$
Note that the elements in each row of the $(p+1)\times(p+1)$
matrix $\mathcal{A}$ are non-negative and sum up to $1$, and thus,
we have a matrix of a finite Markov chain. We need to check that
it is acyclic. Let $\tau(i)=i-1$, and $\sigma(i)=\frac{i}{1-i}$
for $i\in\mathbb{F}_{p}\cup\{\infty\}$. The entry $a_{i,j}^{(m)}$
of $\mathcal{A}^{m}$ is
$$
a_{i,j}^{(m)}=\sum_{i_{1},...,i_{m-1}}a_{i,i_{1}}\cdot a_{i_{1},i_{2}}\cdot...\cdot a_{i_{m-1},j}.
$$
Therefore, we need to check that for some fixed $m$, the
composition of $m$ $\sigma's$ or $\tau's$ leads from any $i$ to
any $j$. One checks directly that for any positive $k$, and
$i,j\in\mathbb{F}_{p}$,
\begin{eqnarray*}
\tau^{p-1-j}\circ\sigma\circ\tau^{k}\circ\sigma\circ\tau^{i-1}(i)&=&j,\\
\tau^{p-1-j}\circ\sigma\circ\tau^{k}(\infty)&=&j,\\
\tau^{k}\circ\sigma\circ\tau^{i-1}(i)&=&\infty;
\end{eqnarray*}
(for $i=0$, we write $\tau^{-1}$ for $\tau^{p-1}$). For each pair $(i,j)$, choose $k$ in order the amount of compositions used to be equal (say, to $m$). Then obviously all entries of $\mathcal{A}^{m}$ are positive, ant this matrix satisfies the conditions of lemma. Since all columns also sum up to $1$, $(\pi_{1},...,\pi_{p+1})$, $\pi_{j}=\frac{1}{p+1}$, $1\leq j\leq p+1$, is the needed eigenvector. This proves the theorem. $\blacksquare$\\

\begin{thm}\label{abel}
Let $\nu\in\mathbb{Z}$ and $z\in\mathbb{Q}_{p}$, and $\text{ord}_{p}(z)<\nu$ (or $z=0$). Then, if $z$ is $p-$adic integer,
$$
\mu(z,\nu)=\frac{1}{p^{\nu}+p^{\nu-1}}.
$$
If $z$ is not integer, $\text{ord}_{p}(z)=-\lambda<0$,
$$
\mu(z,\nu)=\frac{1}{p^{\nu+2\lambda}+p^{\nu+2\lambda-1}}.
$$
For $z=0$, $-\nu\leq0$, we have
$$
\mu(0,-\nu)=1-\frac{1}{p^{\nu+1}+p^{\nu}}.
$$
\end{thm}

\noindent This theorem allows the computation of the associated zeta-function:

\begin{cor} For $s$ in the strip $-1<\Re{s}<1$,
\begin{eqnarray*}
Z_{p}(s)=\int_{u\in\mathbb{Q}_{p}}|u|^{s}d\mu_{p}=
\frac{(p-1)^{2}}{(p-p^{-s})(p-p^{s})},
\end{eqnarray*}
and $Z_{p}(s)=Z_{p}(-s)$.
\end{cor}

\noindent The proof is straightforward. It should be noted that this expression encodes all the values of $\mu(0,\nu)$ for $\nu\in\mathbb{Z}$.
\par\bigskip

\noindent {\bf Proof of Theorem \ref{abel}.} For shortness, when
$p$ is fixed, denote $\text{ord}_{p}(*)$ by $v(*)$. As before, we
want a recurrence relation among the numbers $F_{n}(i,\kappa)$,
$i\in\mathbb{Q}_{+}$. For each integral $\kappa$, we can confine
to the case $i<p^{\kappa}$. If $i=0$, we only consider $\kappa>0$
and call these pairs $(i,\kappa)$ "admissible". We also include
$G_{n}(0,-\kappa)$ for $\kappa\geq 1$, where these values are
defined in the same manner as $F_{n}$, only inversing the
inequality, considering
$\frac{a}{b}\in\mathcal{C}\mathcal{W}^{(n)}$, such that
$v(\frac{a}{b})\leq -\kappa$; the ratio of fractions in the $n$th
generation outside this circle. As before, a fraction
$\frac{a}{b}$ in the $n$th generation generates the fractions
$\frac{a}{a+b}$ and $\frac{a+b}{b}$ in the $(n+1)$st generation.
Let $\tau(i,\kappa)=((i-1)\bmod\,p^{\kappa},\kappa)$. Then for all
admissible pairs $(i,\kappa)$, $i\neq0$, the pair $\tau(i,\kappa)$
is also admissible, and
$$
v(\frac{a+b}{b}-i)=\kappa\Leftrightarrow v(\frac{a}{b}-(i-1))=\kappa.
$$
Second, if $\frac{a}{a+b}=i+p^{\kappa}u$, $i\neq 1$,
$u\in\mathbb{Z}_{p}$, and $(i,\kappa)$ is admissible, then
$$
\frac{a}{b}-\frac{i}{1-i}=\frac{p^{\kappa}u}{(1-i)(1-i-p^{\kappa}u)}.
$$
Since $v(\frac{i}{1-i})=v(i)-v(1-i)$, this is $0$ unless $i$ is an
integer, equals to $v(i)$ if the latter is $>0$ and equals to
$-v(1-i)$ if $v(1-i)>0$. Further, this difference has valuation
$\geq\kappa_{0}=\kappa$, if $i\in\mathbb{Z},i\not\equiv
1\bmod\,p$, valuation $\geq\kappa_{0}=\kappa-2v(1-i)$, if
$i\in\mathbb{Z},i\equiv1\bmod\,p$, and valuation
$\geq\kappa_{0}=\kappa-2v(i)$ if $i$ is not integer. In all three
cases, easy to check, that, if we define
$i_{0}=\frac{i}{1-i}\bmod\,p^{\kappa_{0}}$, the pair
$\sigma(i,\kappa)=^{\text{def}}(i_{0},\kappa_{0})$ is admissible.
For the converse, let $\frac{a}{b}=i_{0}+p^{\kappa_{0}}u$,
$u\in\mathbb{Z}_{p}$. Then
$$
\frac{a}{a+b}-\frac{i_{0}}{1+i_{0}}=\frac{p^{\kappa_{0}}}{(1+i_{0}+p^{\kappa_{0}}u)(1+i_{0})}.
$$
If $i=\frac{i_{0}}{1+i_{0}}$ is a $p-$adic integer, $i\not\equiv 1\bmod\,p$, this has a valuation $\geq\kappa=\kappa_{0}$; if $i$ is a $p-$adic
integer, $i\equiv 1(p)$, this has valuation
$$
\geq\kappa=\kappa_{0}-2v(i_{0})=\kappa_{0}+2v(1-i);
$$
if $i$ is not a $p-$adic integer, this has valuation
$$
\geq\kappa=\kappa_{0}-2v(1+i_{0})=\kappa_{0}+2v(i).
$$
Thus,
$$
v(\frac{a}{a+b}-i)\geq\kappa\Leftrightarrow v(\frac{a}{b}-i_{0})\geq\kappa_{0}.
$$
Let $i=1$. If $\frac{a}{a+b}=1+p^{\kappa}u$, then $\kappa>0$,
$u\in\mathbb{Z}_{p}$, and we obtain
$\frac{a}{b}=-1-\frac{1}{p^{\kappa}u}$,
$v(\frac{a}{b})\leq-\kappa$. Converse is also true. Finally, for
$\kappa\geq1$,
$$
v(\frac{a+b}{b})\leq-\kappa\Leftrightarrow v(\frac{a}{b})\leq-\kappa,
$$
and
$$
v(\frac{a}{a+b})\leq-\kappa\Leftrightarrow v(\frac{a}{b}+1)\geq\kappa.
$$
Therefore, we have the recurrence relations:
\begin{eqnarray}\left\{ \begin{array}{l@{\,}l}
F_{n+1}(i,\kappa)&=\frac{1}{2}F_{n}(\tau(i,\kappa))+\frac{1}{2}F_{n}(\sigma(i,\kappa)),
\text{ if }
(i,\kappa)\text{ is admissible},\\
F_{n+1}(1,\kappa)&=\frac{1}{2}F_{n}(0,\kappa)+\frac{1}{2}G_{n}(0,-\kappa), \kappa\geq1,\\
G_{n+1}(0,-\kappa)&=\frac{1}{2}G_{n}(0,-\kappa)+\frac{1}{2}F_{n}(-1,\kappa), \kappa\geq1.\\
\end{array} \right.\label{pad}
\end{eqnarray}
Thus, we have an infinite matrix $\mathcal{A}$, which is a change
matrix for the Markov chain. If $\mathbf{v}_{n}$ is an infinite
vector-column of $F_{n}'$s and $G_{n}'$s, then
$\mathbf{v}_{n+1}=\mathcal{A}\mathbf{v}_{n}$, and, as before,
$\mathbf{v}_{n}=\mathcal{A}^{n-1}\mathbf{v}_{1}$. It is direct to
check that each column also contains exactly two nonzero entries
$\frac{1}{2}$, or one entry, equal to $1$. In terms of Markov
chains, we need to determine the classes of orbits. Then in proper
rearranging, the matrix $\mathcal{A}$ looks like
\begin{eqnarray*}
\begin{pmatrix}
\mathbf{P}_{1} & 0 & \dots & 0 & \dots\\
0              &\mathbf{P}_{2} & \dots & 0 & \dots\\
\vdots &  & \ddots & \vdots& \vdots\\
0 & 0 & \dots & \mathbf{P}_{s}& 0\\
\vdots & \vdots & \dots & 0 & \ddots
\end{pmatrix},
\end{eqnarray*}
where $\mathbf{P}_{s}$ are finite Markov matrices. Thus, we claim
that the length of each orbit is finite, every orbit has a
representative $G_{*}(0,-\kappa)$, $\kappa\geq 1$, the length of
it is $p^{\kappa}+p^{\kappa-1}$, and the matrix is recurrent (that
is, every two positions communicate). In fact, from the system
above and form the expression of the maps $\tau(i,\kappa)$ and
$\sigma(i,\kappa)$, the direct check shows that the complete list
of the orbit of $G_{*}(0,-\kappa)$ consists of (and each pair of
states are communicating):
\begin{eqnarray*}
&G_{*}(0,-\kappa),\\
&F_{*}(i,\kappa)\quad (i=0,1,2,...,p^{\kappa}-1),\\
&F_{*}(p^{-\lambda}u,\kappa-2\lambda)\quad (\lambda=1,2,...,\kappa-1, u\in\mathbb{N},u\not\equiv 0 \bmod\,p,
u\leq p^{\kappa-\lambda}).
\end{eqnarray*}
In total, we have
$$
1+p^{\kappa}+\sum\limits_{\lambda=1}^{\kappa-1}(p^{\kappa-\lambda}-p^{\kappa-\lambda-1})= p^{\kappa}+p^{\kappa-1}
$$
members in the orbit. Thus, each $\mathbf{P}_{\kappa}$ in the
matrix above is a finite dimensional
$\ell_{\kappa}\times\ell_{\kappa}$ matrix, where
$\ell_{\kappa}=p^{\kappa}+p^{\kappa-1}$. For $\kappa=1$, the
matrix $\mathbf{P}_{1}$ is exactly the matrix of the system
(\ref{padic}). As noted above, the vector column $(1,1,...,1)^{T}$
is the left eigen-vector. As in the previous theorem, it is
straightforward to check that this matrix is irreducible and
acyclic (that is, the entries of $\mathbf{P}_{\kappa}^{n}$ are
strictly positive for sufficiently large $n$). In fact, since by
our observation, each two members in the orbit communicate, and
since we have a move $G_{*}(0,-\kappa)\rightarrow
G_{*}(0,-\kappa)$, the proof of the last statement is immediate:
there exists $n$ such that any position is reachable from another
in exactly $n$ moves, and this can be achieved at the expense of
the move just described. Therefore, all entries of
$\mathbf{P}_{\kappa}^{n}$ are strictly positive. Thus, the claim
of the theorem follows from the lemma above. $\blacksquare$

\par\bigskip

\noindent {\bf Acknowledgements.} The second author thanks J\"urgen Sander and Jan-Hendrik de Wiljes for their interest and valuable remarks.

\par\bigskip

\noindent Giedrius Alkauskas, Department of Mathematics and Informatics, Vilnius University, Naugarduko 24, 03225 Vilnius, Lithuania\\
 giedrius.alkauskas@maths.nottingham.ac.uk
\smallskip

\noindent
J\"orn Steuding, Department of Mathematics, W\"urzburg University, Am Hubland, 97\,218 W\"urzburg, Germany\\
steuding@mathematik.uni-wuerzburg.de
\smallskip

\end{document}